\documentclass{amsart}
\usepackage{amscd}
\usepackage{amsmath}
\usepackage{amsxtra}
\usepackage{amsfonts}
\usepackage{amssymb}

\newtheorem{theorem}{Theorem}[section]
\newtheorem{corollary}[theorem]{Corollary}
\newtheorem{lemma}[theorem]{Lemma}
\newtheorem{proposition}[theorem]{Proposition}
\theoremstyle{definition}
\newtheorem{definition}[theorem]{Definition}
\newtheorem{remark}[theorem]{Remark}

\newtheorem{example}[theorem]{Example}
\theoremstyle{remark}

\renewcommand{\theclaim}{\textup{\theclaim}}

\newtheorem*{acknowledgements}{Acknowledgements}

\numberwithin{equation}{section}

\def\openone

{\mathchoice

{\hbox{\upshape \small1\kern-3.3pt\normalsize1}}

{\hbox{\upshape \small1\kern-3.3pt\normalsize1}}

{\hbox{\upshape \tiny1\kern-2.3pt\SMALL1}}

{\hbox{\upshape \Tiny1\kern-2pt\tiny1}}}

\makeatletter

\newbox\ipbox

\newcommand{\diracb}[1]{\left\langle #1\mathrel{\mathchoice

{\setbox\ipbox=\hbox{$\displaystyle \left\langle\mathstrut
#1\right.$}

\vrule height\ht\ipbox width0.25pt depth\dp\ipbox}

{\setbox\ipbox=\hbox{$\textstyle \left\langle\mathstrut
#1\right.$}

\vrule height\ht\ipbox width0.25pt depth\dp\ipbox}

{\setbox\ipbox=\hbox{$\scriptstyle \left\langle\mathstrut
#1\right.$}

\vrule height\ht\ipbox width0.25pt depth\dp\ipbox}

{\setbox\ipbox=\hbox{$\scriptscriptstyle \left\langle\mathstrut
#1\right.$}

\vrule height\ht\ipbox width0.25pt depth\dp\ipbox}

}\right. }

\newcommand{\dirack}[1]{\left. \mathrel{\mathchoice

{\setbox\ipbox=\hbox{$\displaystyle \left.\mathstrut
#1\right\rangle$}

\vrule height\ht\ipbox width0.25pt depth\dp\ipbox}

{\setbox\ipbox=\hbox{$\textstyle \left.\mathstrut
#1\right\rangle$}

\vrule height\ht\ipbox width0.25pt depth\dp\ipbox}

{\setbox\ipbox=\hbox{$\scriptstyle \left.\mathstrut
#1\right\rangle$}

\vrule height\ht\ipbox width0.25pt depth\dp\ipbox}

{\setbox\ipbox=\hbox{$\scriptscriptstyle \left.\mathstrut
#1\right\rangle$}

\vrule height\ht\ipbox width0.25pt depth\dp\ipbox}

} #1\right\rangle}

\newcommand{\bz}{\mathbb{Z}}
\newcommand{\br}{\mathbb{R}}
\newcommand{\bc}{\mathbb{C}}

\newcommand{\bn}{\mathbb{N}}

\usepackage{graphicx}

\def\blfootnote{\xdef\@thefnmark{}\@footnotetext}

\begin{document}
\title[Harmonic analysis]{Harmonic analysis and dynamics for affine iterated function systems}
\author{Dorin Ervin Dutkay and Palle E.T. Jorgensen}
\blfootnote{Work supported in part by the National Science Foundation.}
\address{Department of Mathematics\\
The University of Iowa\\
14 MacLean Hall\\
Iowa City, IA 52242-1419\\
U.S.A.\\} \email[Dorin Ervin Dutkay]{ddutkay@math.rutgers.edu}
\email[Palle E.T. Jorgensen]{jorgen@math.uiowa.edu}
\begin{abstract}
We introduce a harmonic analysis for a class of affine iteration models in $\br^d$. Using Hilbert-space geometry, we develop a new duality notion for affine and contractive iterated function systems (IFSs) and we construct some identities for the Fourier transform of the measure corresponding to infinite Bernoulli convolutions.
\end{abstract}
\maketitle
\tableofcontents

\section{Introduction}
There has been a recent increased interest in basis constructions outside
the traditional setting of harmonic analysis. The setting which so far has proved
more amenable to an explicit analysis with basis functions involves a mix of
analysis and dynamics, and it typically goes beyond the standard and more familiar
setting of orthonormal bases consisting of Fourier frequencies. The context of
frames in Hilbert space (see, e.g., \cite{ALTW04, BoPa05, BPS03}) is a case in point.

          As is well known, the classical setting of Fourier presupposes an
underlying group and Haar measure. But in dynamics, this is not available, and
there is a variety of interesting cases which go beyond the group setting, and
which nonetheless are amenable to harmonic basis constructions; see, e.g., \cite{JoPe98,
DuJo03,DuJo05}. These constructions are recursive in nature.  At the same
time, these recursive models arise in applications exhibiting a suitable
{\em scale-similarity}. Their consideration combines classical ideas from infinite
convolution with ideas from dynamics of a more recent vintage.

         In the simplest and still useful non-group case, the similarity may be
defined by a scaling of vectors in $\br^d$ with a fixed $d$ by $d$ matrix. Our sets and
measures then arise as attractors of an associated dynamical system. The sets and
measures in our present non-group context arise from a certain algorithmic
construction on a fixed finite system of affine mappings in $\br^d$. Consideration of bases in this context involves a new
viewpoint and a combination of ideas from dynamics and harmonic analysis.

Our aim in this paper is to develop a notion of cyclicity from dynamical systems and operator theory, and to apply it to
geometric iteration systems. Specifically, we create a harmonic analysis on the class of affine iterated function systems
(IFS). Our present harmonic analysis is based on Hilbert-space geometry and a certain duality notion, but it is more
versatile than the earlier such notions from fractal Fourier analysis; see especially \cite{JoPe98}. As is well known,
IFSs provide us with prototypes of fractals, and they have for some time served as successful models for symbolic dynamics,
e.g., \cite{BKS91}. More generally, fractals of the affine IFS-type arise in several areas in symbolic dynamics;
see, e.g., \cite{Edg98}.

Specifically, we develop in this paper a concept of cyclicity which generalizes a harmonic analysis which was studied in
earlier papers by the second named author, S. Pedersen, and R. Strichartz, among others. While it was
shown in \cite{JoPe98} that certain affine and contractive iterated function systems in $\br^d$ admit Hilbert-%
space-orthonormal bases (ONBs), it turns out that the conditions for existence of ONBs in affine IFSs is somewhat
restrictive. The affine systems under discussion include Cantor sets of Sierpinski type, realized in $\br^d$.
\par
In Section \ref{background} below, we recall how a fixed iterated function system (IFS) $(\tau_i)_i$ gives rise to a unique attractor $X$. As a result, we get induced maps in $X$, also denoted by $\tau_i$. We further recall how the index labels for the system $(\tau_i)$ introduce an encoding space $\Omega$ for the points in $X$. Using the natural shift map of $\Omega$ and the notion of periodicity, we may then introduce associated cycles $C$ as subsets of $X$. A certain weight function $W$ is fixed and defined on $X$; and it allows us to introduce transition probabilities associated with $(\tau_i)$, a transition operator $R_W$, and a path-space measure $P_x = P_x^{(W)}$. We further use $W$ to define certain $W$-cycles which are key to understanding the solutions to the equation $R_W h = h$. These solutions $h$ are called harmonic functions, or $R_W$-harmonic functions.
\par
   Our main results Theorems \ref{proph-0=1}, \ref{thunique} are about the harmonic functions for the special IFSs of Bernoulli type. Section \ref{reps} deals with some geometric issues used in Section \ref{harm}. Section \ref{meas} is about general IFSs. 
\par
It turns out that the introduction of a harmonic analysis for these systems adds to their use, and it shows that
seemingly similar Cantor sets can have drastically different harmonic analysis: For example, the usual Cantor subset
of the line, arising from subdivision by 3, is vastly different from the corresponding Cantor construction based instead
on subdivision by 4; see \cite{JoPe98} for details. In this paper, we identify a notion of extremality for certain finite
cycles in IFSs. Among these cycles, we identify a small subclass of cycles where the given weight function $W$ attains
its maximum. The choice of weight function $W$ depends on and is intimately tied to the geometry of the IFS under
discussion. In fact, $W$ serves to create a duality of affine IFSs via a certain complex Hadamard matrix; see, e.g.,
\cite{JoPe96}.
 When the IFS and $W$ are specified, we call the corresponding extreme cycles $W$-cycles (Definition \ref{defcycle}). We point out
with examples the significance of these $W$-cycles; and we find all the $W$-cycles for particular affine fractal systems.
It turns out that the case when an affine IFS has the ONB property (such as the Cantor set arising by quarter partitions
and skipping every second subinterval) corresponds to the situation in our present setup when the only extreme $W$-cycle
is a one-cycle. This last fact holds more generally; but as we show, not all systems for which the only $W$-cycle is a
one-cycle admit an ONB of Fourier type.

\par
A $W$-cycle is a cycle $C$ for $(X, \tau_i)$ such that $W(y) = 1$ for all $y$ in $C$ (Definition \ref{defcycle}). The function $W$ determines a transition operator $R_W$ (Definition \ref{defruelle}) and 
a random walk on a space of paths beginning at points in $x$. Hence, each $x$ induces a Kolmogorov measure $P_x$ on $\Omega$ (equation (\ref{eq3.diez})).
We show that each $P_x$ is countably atomic, i.e., is supported on a countable set $S$ in $\Omega$, and each point in
$S$ is an atom for $P_x$ (Corollary \ref{cor2-1}, Theorem \ref{proph-0=1}). In particular, we give the precise structure of this $P_x$ for a classical one-parameter family of infinite Bernoulli convolution measures $\nu_\lambda$ (Examples \ref{exbern}, \ref{exbern2}), and associated support sets $X_\lambda$. Using this, we further show that several properties of $\nu_\lambda$ are determined directly
from the pair  $(P_x , R_W)$. For example, we find the support $X_\lambda$ of $\nu_\lambda$ (Proposition \ref{lembern1}); and we show that for all $\lambda$, the absolute square of the Fourier transform
of $\nu_\lambda$ is the unique solution to a functional equation defined directly from $P_x$ (Theorem \ref{thunique}). The discrete set $S$ which
supports the measures $P_x$ encodes the union of the $W$-cycles (Lemma \ref{lemnozerois1}); and our proofs are based on a detailed analysis of
the $W$-cycles for the duality of the Bernoulli systems. In particular, we find all the $W$-cycles for the Bernoulli
systems and we show (Theorem \ref{proph-0=1}) that there is a one-to-one correspondence between the $W$-cycles $C$ and certain extreme solutions
$h(x)=h_C(x)=P_x(\{C\mbox{-cyclic paths}\})$ to the eigenvalue problem $R_W h = h$, called harmonic functions (see Definition \ref{defNC}). We show that there is a specific set $D$ of rational values of
$\lambda$, such that if $\lambda$ is in $(0,1) \setminus D$, then $H_\lambda:=\{ h \in C(X_\lambda)\, |\,  R_W h = h\}$ is
one-dimensional, and if $\lambda$ is in $D$, then $\dim H_\lambda = 2$.  For all $\lambda$, we show that
$\sum_C h_C = 1$ on $X_\lambda$.

When $W$ is assumed normalized, we obtained in \cite{DuJo05} that, under certain conditions, $\sum_Ch_C=1$, i.e., the constant function
``one''. But we also gave examples in $\br^d$ when this does not hold. As an application of our $W$-cycles, we show that for the Bernoulli IFSs,
this sum {\it is} the constant function 'one' for {\it all} $\lambda$ in $(0, 1)$.
Here we address this issue and related questions for the harmonic analysis of Bernoulli IFSs, and our discussion
is divided up into the two cases: (i) $\lambda$ in $(0,\frac12)$, the fractal case, and (ii) $\lambda\geq\frac12$.
It is the second case that holds more surprises.
\par
       Our first result (see Theorem \ref{lemnocycles}) for these
Bernoulli systems is that when $\lambda$ is in $[1/2, 1)$ there are
no $W$-cycles of minimal period $p > 1$.  If $\lambda > 1/2$,  and $\lambda\not\in D:=\{ 1-\frac1{2n}\,|\,n\in\bn\}$, then
the only $W$-cycle is $C = \{0\}$ as a subset of the interval
$X_\lambda$ which supports the system.  The corresponding
eigenfunction, or harmonic function, is then   $h_0 (x) =
P_x(\mathbf{N}_0)$, and $\mathbf{N}_0$ is the set of paths that end in an infinite repetition of the trivial cycle $\{0\}$. So either $h_0$ is the constant function $1$, or
 $\dim H_\lambda > 1$.

In Theorem \ref{proph-0=1}, we show that if $\lambda  > \frac12$
falls in the complement of the countable and exceptional set
$D$, then we may conclude
that $h_0$ has no zeroes.  Then by a theorem from \cite{DuJo05},
$h_0$ must be the constant function 'one', and indeed we get
$\dim H_\lambda = 1$. The significance of this may be understood from our result,
Theorem \ref{thunique}. This theorem shows that $|\hat\nu_\lambda|^2$ is the unique
solution to a certain functional equation, i.e., the modulus-square of
the Fourier transform of the Bernoulli measure $\nu_\lambda$ is determined
this way for all $\lambda$ in the complement of $D$.

For $\lambda$ in the exceptional set $D$, the space $H_\lambda$ is also understood (see part two in Theorem \ref{proph-0=1}):
In that case, $h_0$ is not a constant function; the space $H_\lambda$ is then $2$-dimensional, and it
is generated by two distinct $W_\lambda$-cycles $C_0$ and $C_1$. Each of the two cycles defines a harmonic function,
$h_0$ and $h_1$, respectively. The two functions form a basis for $H_\lambda$, and $h_0  +  h_1  = 1$.

\section{Definitions}\label{background}
\begin{definition}\label{def4-1}
An affine IFS in $\br^d$ is determined by
\begin{itemize}
\item
$R$: a $d\times d$ matrix such that the eigenvalues $\lambda$ satisfy $|\lambda|>1$;
\item
$B\subset\br^d$ a finite subset.
\end{itemize}
Then we define
\begin{equation}\label{eqtau1}
\tau_b(x):=R^{-1}(x+b),\quad(x\in\br^d),
\end{equation}
and note that there is a unique compact subset $X_B\subset\br^d$ such that
\begin{equation}\label{eqXB2}
X_B=\bigcup_{b\in B}\tau_b(X_B).
\end{equation}
\end{definition}
In fact $X_B$ consists of the points
\begin{equation}\label{eqexpa3}
x=\sum_{j=1}^\infty R^{-j}b_j,\quad\mbox{as }b=(b_1,b_2,\dots )\in\prod_{\bn}B\quad(=:B^{\bn}).
\end{equation}
Set $\Omega=\Omega_B=B^{\bn}$;
\begin{itemize}
\item $\mathcal{W}:=$all finite words in $B$;
\item $\mathcal{W}_k:=\{(b_1\dots b_k)\,|\,b_i\in B\}$;
\item If $b\in\mathcal{W}$, let $|b|$ denote the length of $b$, i.e., $|(b_1\dots b_k)|=k$.
\end{itemize}
For $b\in\Omega_B$, set
\begin{equation}\label{eqpi4}
\pi_R(b):=\sum_{j=1}^\infty R^{-j}b_j.
\end{equation}
For $k\in\bn$, and $b\in\mathcal{W}_k$, set
\begin{equation}\label{eqpi5}
\pi_R(b):=\sum_{j=1}^kR^{-j}b_j.
\end{equation}

 When an affine IFS $(X,\tau_i)$ is given, our duality is then determined by a certain
function $W \colon  X\rightarrow [0,\infty)$ satisfying
\begin{equation}\label{eqnorm}
\sum_{i=1}^NW(\tau_ix)=1,
\end{equation}
and by an associated transfer operator $R_W$ acting on functions on $X$ (Definition \ref{defruelle}). We show that $R_W$ has properties analogous
to those of finite Perron-Frobenius non-negative matrices.
\par
For a non-negative function $W$ on $X$ satisfying (\ref{eqnorm})
following Kolmogorov, one can define probability measures $P_x$ on
$\Omega$, $x\in X$, such that, for a function $f$ on $\Omega$ which
depends only on the first $n$ coordinates,
\begin{multline}\label{eq3.diez}
\int_\Omega f(\omega)\,dP_x(\omega)=P_x(f)\\
=\sum_{\omega_1,\dots ,\omega_n}W(\tau_{\omega_1}x)W(\tau_{\omega_2}\tau_{\omega_1}x)\cdots W(\tau_{\omega_n}\cdots\tau_{\omega_1}x)f(\omega_1,\dots ,\omega_n).
\end{multline}
We will assume throughout that the {\it normalization} (\ref{eqnorm}) holds. In particular $W_\lambda$ in (\ref{eqWB14}) for the $\lambda$-Bernoulli system satisfies (\ref{eqnorm}).

      The main question addressed in this paper arose in, or is motivated by, a number of earlier investigations on iterated function
systems (IFS); see, e.g., \cite{Edg98}, \cite{JoPe96}, \cite{JoPe98}, \cite{So95}, \cite{PeSo96}, \cite{KSS03}, and \cite{Jor04}. For additional details on the measures $P_x$ occurring in (\ref{eq3.diez}), see, for
example, \cite{Jor06}.
The Bernoulli convolutions (see details below) are of interest in their own right, and serve as models
for general structures in dynamics. They were studied in the thirties by P. Erd\"os \cite{Er39}
and subsequently by many others.

\par
Infinite Bernoulli convolutions (Example \ref{exbern}, Proposition \ref{lembern1}) form a subclass of contractive iterated function systems (IFS) of affine type (Definition \ref{def4-1}).
Each such affine IFS $(X, (\tau_i)_{i=1,\dots ,N})$ has a unique measure $\nu_p$ on $X$ which assigns prescribed
probabilities $p_i$ to the branch maps $\tau_i$. We establish a duality for these measures, and we show in particular
that each $\nu_p$ is non-atomic (Corollary \ref{cornoatombern}).
\par
Let $\lambda$ be in the open interval $(0, 1)$. Let $N=2$, $p_i=\frac12$,  let the two $\tau$-maps be $\tau_{\pm}(x):=\lambda x\pm 1$,
and denote the corresponding measure by $\nu_\lambda$ (Example \ref{exbern2}). This measure $\nu_\lambda$ is the distribution of the random
variable $x(\omega)=\sum_{k=1}^\infty\pm\lambda^k$ with independent and random assignment of $\pm$ signs.
Moreover, $\nu_\lambda$ is an infinite convolution measure, studied first by P. Erd\"os, and is called a Bernoulli
convolution. Let $\Omega$ be the space of all configurations $\omega$ of signs in $x$.
We show that each Bernoulli system is {\it one side} in a duality of affine IFSs, and that this duality holds more generally
for contractive affine IFSs.

The IFS $\tau_\pm(x)=\lambda x\pm1$ defined for fixed $\lambda$ is called a {\it Bernoulli IFS\/}; and for each $\lambda$ in $(0, 1)$, there is a natural
choice for the function $W =W_\lambda$, see equation (\ref{eqWB14}), and an associated process $P_x$.

\begin{definition}\label{defruelle}
A third essential ingredient for the process is a transfer operator $R_W$ depending on $W$. For functions $f$ on $X$ we set
$$(R_Wf)(x)=\sum_{i=1}^NW(\tau_ix)f(\tau_ix)\quad(x\in X).$$
We denote the operator $R_W$, the {\it Ruelle operator}, or a transfer operator, to stress the dependency on the function $W$.
Each $R_W$ defines special
{\it harmonic functions}  $$H = \{h \in C(X)\, |\, R_W h = h \}.$$
\end{definition}

The main results in this paper concern the real part of the peripheral spectrum of this operator $R_W$; see especially Theorem \ref{proph-0=1} below. But to state and prove our spectral results for $R_W$, we must first introduce a certain encoding map for our IFSs, and an associated family of finite cycles, called $W$-cycles.
\begin{definition}\label{defcycle}
{\bf[Cycles and $W$-cycles]}\enspace
If a point $x\in X$ has the property that
$\tau_{\omega_{p}}\cdots \tau_{\omega_1}x=x$ for some
$\omega_1,\dots ,\omega_p\in\{1,\dots ,N\}$, then the set
$$C:=\{x,\tau_{\omega_1}x, \dots  , \tau_{\omega_{p-1}}\cdots \tau_{\omega_1}x\}$$
is called a {\it cycle}. If $p$ is minimal with this property, then $p$ is
called the length of the cycle and $C$ is called a $p$-cycle. We can identify the cycle $C$ with the finite word
$\omega_1\dots \omega_p$.

Every cycle with period $p$ is also periodic with a period that is a
multiple of $p$. We say that a cycle $C$ has minimal period $p$ if
$C$ does not close up with a period smaller than $p$, i.e., with a
period that is a divisor in $p$. Unless specified otherwise, we will
assume that the term $p$-cycle refers to the minimal period.\par
When $W$ is given, we are interested in those cycles $C$ satisfying
$W(x) = 1$ for all $x$ in $C$. These cycles are called {\it $W$-cycles}. And we
will use the term $W$-$p$-cycle for a $p$-cycle which is also a
$W$-cycle.
\end{definition}

\begin{definition}\label{defencoding}
Let $\Omega=\{1,\dots ,N\}^{\bn}$. We say that $\Omega$ is an encoding of $(\tau_i)_{i=1}^N$ if for all $\omega=(\omega_1\omega_2\dots )\in\Omega$
\begin{equation}\label{eqdef33}
\bigcap_{n=1}^\infty\tau_{\omega_1}\tau_{\omega_2}\cdots \tau_{\omega_n}(X)
\end{equation}
is a singleton $\pi(\omega)$. (This condition is
automatic if the maps $\tau_i$ are strictly contractive, but it is known
to hold also for important classes of IFS which are not contractive.)
\par
Suppose an IFS has an encoding $\pi$. Then $\pi\colon \Omega\rightarrow X$ maps onto the attractor of the IFS $X_\tau$, and $\pi$ is called the {\it encoding mapping}.
\end{definition}

If $(\omega_1\dots \omega_p)$ is a word in the alphabet $\{1,\dots ,N\}$, then
$$(\omega_1\dots \omega_p)^\infty=(\underline{\omega_1\dots \omega_p}\,\underline{\omega_1\dots \omega_p}\dots )$$
is the point in $\Omega$ resulting from an infinite repetition of $(\omega_1\dots \omega_p)$. It generates a $p$-cycle under the shift $\sigma$ in $\Omega$, i.e.,
\begin{equation}\label{eqshift}
\sigma(\eta_1\eta_2\eta_3\dots ):=(\eta_2\eta_3\eta_4\dots ),
\end{equation}
i.e.,
$$\sigma^p((\omega_1\dots \omega_p)^\infty)=(\omega_1\dots \omega_p)^\infty.$$
This $p$-cycle will be denoted $C(\omega_1\dots \omega_p)$, and it is said to be generated by the word $(\omega_1\dots \omega_p)$.

\begin{definition}\label{defNC}
For each $W$-cycle $C$, consider the subset $\mathbf{N}_C$ in $\Omega$
consisting of paths that end up in $C$ after a finite number of
backwards shifts, i.e., $\mathbf{N}_C$ consists of all finite words, followed by an infinite repetition of the fixed word representing $C$. That is, if $C$ is generated by the word $(k_1\dots k_p)$, then
$$
\mathbf{N}_C:=\{\omega_1\dots \omega_n(k_1\dots k_p)^\infty\,|\,\omega_1,\dots ,\omega_n\in\{1,\dots ,N\},n\in\bn\}.
$$
We define the function
\begin{equation}\label{eqhc}
h_C(x) = P_x(\mathbf{N}_C).
\end{equation}
\end{definition}
 Then $h_C$ is harmonic and continuous (see \cite[Proposition 5.8, Remark 5.9]{DuJo05}), and the question is if the space $H$ of continuous harmonic functions is generated by the $W$-cycles.\par
 Note the following
distinct parts of the conclusion: First that the function $h_C(x)=P_x(\mathbf{N}_C)$ is continuous, and secondly,
that it satisfies the equation $R_W h_C = h_C$. The reader is also referred to \cite{DuJo05} for the complete arguments of this part.
It turns out that properties for the functions $h_C$, such as existence of zeroes,
can be checked; and that these properties in turn have implications both for the measures $P_x$
and for the strongly invariant IFS measures $\nu$.

\section{$W$-cycles for infinite Bernoulli convolutions}\label{seccycles}

\begin{example}\label{exbern}
Bernoulli's infinite convolution; see \cite{SSU01}, \cite{So95},
\cite{PeSo96}, \cite{PeSo00}, and \cite{Er39}. Let $d=1$. Let $\lambda\in(0,1)$, and set $R=\lambda^{-1}$, $B=\{0,1\}$. Then $\Omega_B$ is the usual
probability space $\Omega_B=\{0,1\}^{\bn}$.
\end{example}

\begin{example}\label{exbern2}
A dual pair of Bernoulli systems.
\par
Let $\lambda\in(0,1)$ and set $R=\lambda^{-1}$. We consider IFSs corresponding to two sets
\begin{equation}\label{eqblambda}
B(\lambda):=\{\pm\lambda^{-1}\},\quad\mbox{and}\quad L(\lambda):=\left\{0,\frac14\right\},
\end{equation}
i.e.,
\begin{equation}\label{eqtaupm}
\tau_{\pm}(x):=\lambda x\pm 1;
\end{equation}
and
\begin{equation}\label{eqtau01/4}
\tau_0(x):=\lambda x,\quad \tau_{1/4}(x):=\lambda\left(x+\frac14\right).
\end{equation}
The choice of $B$, and the term {\it duality}, are motivated by the general theory of affine IFSs as developed in \cite{JoPe96}, \cite{JoPe98}, and \cite{DuJo05}, as well as in other recent papers on the harmonic analysis of IFSs.
\par
Setting
\begin{equation}\label{eqWB14}
W(x):=W_B(x)=\cos^2(2\pi\lambda^{-1}x),
\end{equation}
and
\begin{equation}\label{eqRW15}
(R_Wf)(x):=W(\tau_0x)f(\tau_0x)+W(\tau_{1/4}x)f(\tau_{1/4}x),
\end{equation}
we see that
\begin{equation}\label{eqRW16}
R_Wf(x)=\cos^2(2\pi x)f(\lambda x)+\sin^2(2\pi x)f(\lambda(x+1/4)),
\end{equation}
and
\begin{equation}\label{eqRW17}
R_W\mathbf1=\mathbf1,
\end{equation}
where $\mathbf1$ denotes the constant function ``one''.
\par
We will see in Proposition \ref{lembern1} that if $\lambda\geq1/2$, then the two attractors $X_{B(\lambda)}$ and $X_{L(\lambda)}$ consist of the following two intervals:
$$X_{B(\lambda)}=\left[\frac{-1}{1-\lambda},\frac{1}{1-\lambda}\right]\quad\mbox{and}\quad X_{L(\lambda)}=\left[0,\frac{\lambda}{4(1-\lambda)}\right].$$
\end{example}
\par
The classical theory of Bernoulli convolutions \cite{PSS00} has a point of contact with the scaling equation
\cite{Dau92} from wavelet theory. To illustrate this point, consider functions $f$ on $\br$, and the equation
\begin{equation}\label{eqinsdz1}
f(x) =\frac{1}{2\lambda}\left(f\left(\frac{x-1}{\lambda}\right) + f\left(\frac{x+1}{\lambda}\right)^{\mathstrut}\right).
\end{equation}
Its integrated form is as follows:
\begin{equation}\label{eqinsdz2}
F(x) =\frac12\left(F\left(\frac{x-1}{\lambda}\right)+ F\left(\frac{x+1}{\lambda}\right)^{\mathstrut}\right).
\end{equation}
\par
Then, by elementary multi-scale theory, the special case of $\lambda =\frac12$ in (\ref{eqinsdz1}) coincides
with Haar's scaling identity; see, e.g., \cite{Dau92}. To see this, recall that Haar's scaling function is $f =\frac14\chi_{[ -2, 2]}$.
 A calculation shows that  (\ref{eqinsdz2}) is solved by
$$F(x)=\left\{\begin{array}{ccc}
0&\mbox{ if }& x < -2,\\
\frac14(x+2)&\mbox{ if }&-2\leq x < 2,\\
1&\mbox{ if }&2\leq x.
\end{array}\right.$$
In general, if $\lambda$ is in $(0, 1)$ and $\nu_\lambda$ denotes the Bernoulli measure, then the cumulative
distribution function $F_\lambda(x) : =\nu_\lambda((-\infty, x])$ solves (\ref{eqinsdz2}), and the special case
$\lambda=\frac12$ coincides with $F$ above.

\par
Let
$$W(x)=\cos^2\left(\frac{2\pi x}{\lambda}\right)\quad(x\in\br).$$
We saw in Example \ref{exbern2} that, for the IFS $\tau_0(x)=\lambda x$, $\tau_{1/4}(x)=\lambda(x+\frac14)$, the weight function
$W$ satisfies the normalization condition
$$R_W1=1.$$
Also the attractor of this IFS is $X_{L(\lambda)}=[0,\frac{\lambda}{4(1-\lambda)}]$ if $\lambda>1/2$. We are now ready to list the $W$-cycles for the Bernoulli convolutions for $\lambda$ in the interval $(0,1)$.

We will show that for the Bernoulli IFS systems with parameter
$\lambda$, there are no $W$-$p$-cycles for $p > 1$, when $\lambda$ is
in the open interval $(0, 1)$.

\begin{theorem}\label{lemnocycles}
Suppose $\lambda\in(0,1)$. Then there are no $W$-$p$-cycles for $p>1$.
\end{theorem}

\begin{proof}We relabel $\tau_1:=\tau_{1/4}$.
Assume a priori that $C$ is some $W$-$p$-cycle.
Since it is not a $W$-$1$-cycle, the word $\omega$ which generates
$C$ must contain both the letters $0$ and $1$ from the alphabet
$\{0,1\}$.

First note that
\begin{equation}\label{eqW=1}
W(x)=1\mbox{ only if }2x/\lambda\in\bz.
\end{equation}
\par
Let $x$ be now a point in a $W$-$p$-cycle, with $p>1$.
A $p$-cycle is generated by a word $\omega:=(\omega_1\dots \omega_p)$ of length $p$ over the alphabet $\{0,1\}$, i.e., $x=\pi_\lambda(\omega^\infty)$, where
$$\omega^\infty=\underbrace{(\omega\,\omega\,\omega\dots\dots)}_{\mbox{infinite repetition}}.$$
But
$$\pi_\lambda(\omega^\infty)=\frac{1}{4}\lambda(\omega_1+\omega_2\lambda+\dots +\omega_p\lambda^{p-1})\frac{1}{1-\lambda^p}.$$
Therefore, with equation (\ref{eqW=1}), we have that
\begin{equation}\label{eqnocyc1}
(\omega_1+\omega_2\lambda+\dots +\omega_p\lambda^{p-1})\frac{1}{2(1-\lambda^p)}\in\bz.
\end{equation}
\par
We claim that $\lambda$ must be rational. Indeed, after a cyclic permutation we may assume $\omega_1=0$. Then look at the next point in the cycle, which is generated by the word $(\omega_2\omega_3\dots \omega_p\omega_1)$. Clearly, this point will be $\lambda x$. Using the fact that $x$ and $\lambda x$ are points in a $W$-cycle, so $W(x)=W(\lambda x)=1$, we obtain that
$$2x/\lambda\in\bz,\quad\mbox{and }2x\in\bz.$$
This implies that $\lambda$ is rational.
\par
Let $\lambda=a/b$ with $a,b\in\bn$, $a<b$ and $a,b$ mutually prime. Using (\ref{eqnocyc1}), we obtain that
$$b(\omega_1b^{p-1}+\omega_2b^{p-2}a+\dots +\omega_pa^{p-1})\frac{1}{2(b^p-a^p)}\in\bz.$$
If a prime number divides $b^p-a^p$ then it cannot divide $b$, because, if it does, it must also divide $b^p-(b^p-a^p)=a^p$. Then it would divide both $a$ and $b$, which contradicts the fact that $a$ and $b$ are mutually prime.
\par
Therefore we must have that $b^p-a^p$ divides $\omega_1b^{p-1}+\omega_2b^{p-2}a+\cdots \omega_{p}a^{p-1}$. However
$$\omega_1b^{p-1}+\omega_2b^{p-2}a+\cdots \omega_{p}a^{p-1}<b^{p-1}+b^{p-2}a+\dots +a^{p-1}=\frac{b^p-a^p}{b-a}\leq b^p-a^p.$$
The inequality is strict, because at least one of the $\omega_i$'s is $0$. We have reached a contradiction, and the proof is complete. We conclude that when $\lambda\in(0,1)$, there are
no $W$-cycles of (minimal) period bigger than $1$.
\end{proof}

\section{A harmonic function}\label{harm}
In this section, we address the problem of finding all the harmonic
functions for the infinite Bernoulli convolutions defined for
$\lambda$ in the open interval $(0, 1)$. For every $\lambda$ in $(0,
1)$, we are studying the corresponding Bernoulli-IFS. The present section has two main results; first a lemma
(Lemma \ref{lemnozerois1}) about the zeroes of IFS-harmonic
functions; and secondly our main result (Theorem \ref{proph-0=1}):
If $\lambda$ is assumed irrational, then the constant functions are
the only continuous IFS-harmonic functions. Specifically, the
eigenspace $H_\lambda$ in $C(X_\lambda)$ for the transfer operator
$R_{W_\lambda}$ with eigenvalue $1$ is spanned by the constant
function $h = 1$.

         Recall that $R_{W_\lambda}$ denotes the basic Ruelle operator
         (\ref{eqRW16}), and $W=W_\lambda$ is defined in (\ref{eqWB14}).
         And note further that we are addressing only the question of continuous solutions to
         $R_{W_\lambda} h = h$. In fact, it is intriguing that the case when $\lambda$ is rational is the one
         that is more subtle as far as the IFS-harmonic functions are concerned. Finally, if measurable eigenfunctions
         are admitted into the analysis, then there are many more solutions; see, e.g., \cite{Ba00}, \cite{Jor06}, \cite{JoPe98},
         \cite{Nu01}, and \cite{LMW96}. We shall now apply the encoding $\mathbf{N}_C$ from Definition \ref{defNC} of the
possible $W$-cycles to the Bernoulli systems. We already saw that for
the Bernoulli systems, the only possibilities for $W$-$p$-cycles are $p=1$;
and that then only $C_0$ and $C_1$ are possible. Which combinations
actually occur depends on a partition of the value of $\lambda$ (details
below in Lemma \ref{lemnozerois1}).
\par

Let
$$\mathbf{N}_0:=\{\omega=(\omega_1\omega_2\dots )\in\{0,1/4\}^\bn\,|\,\omega_n=0,\mbox{
for }n\mbox{ big enough}\},$$ i.e., $\mathbf{N}_0$ is the set of
infinite paths that end in $000\dots $.
\par
We saw that if $\lambda=1-\frac{1}{2n}$ for some $n\in\bn$, then there is another $W$-cycle of length $1$,
namely $\{\frac{\lambda}{4(1-\lambda)}\}$.
In this case we also define
$$\mathbf{N}_1:=\{\omega=(\omega_1\omega_2\dots )\in\{0,1/4\}^\bn\,|\,\omega_n=1/4,\mbox{
for }n\mbox{ big enough}\},$$ i.e., $\mathbf{N}_1$ is the set of
infinite paths that end in $\frac14\frac14\frac14\dots $.
\par
Then let, for $x\in X_L$,
\begin{equation}\label{eqh-0}
h_0(x):=P_x(\mathbf{N}_0)=\sum_{\omega\in\mathbf{N}_0}P_x(\{\omega\})=\sum_{\omega\in\mathbf{N}_0}\prod_{k=1}^\infty W(\tau_{\omega_k}\cdots \tau_{\omega_1}x).
\end{equation}
In the case when $\lambda=1-\frac1{2n}$, we also define
\begin{equation}\label{eqh-1}
h_1(x):=P_x(\mathbf{N}_1)=\sum_{\omega\in\mathbf{N}_1}P_x(\{\omega\})=\sum_{\omega\in\mathbf{N}_1}\prod_{k=1}^\infty W(\tau_{\omega_k}\cdots \tau_{\omega_1}x).
\end{equation}

\begin{lemma}\label{lemnozerois1}
Consider the functions $h_0$ and $h_1$ defined above.
\begin{enumerate}
\item
The function $h_0$ is continuous and harmonic for the transfer
operator $R_W$, i.e., $R_Wh_0=h_0$. If $h_0$ has no zeroes then
$h_0=1$.
\item
If $\lambda=1-\frac1{2n}$ then $h_1$ is continuous and
harmonic for the transfer operator $R_W$. Both functions then have zeroes, and both are non-constant. If $h_0+h_1$ has no zeroes
then $h_0+h_1=1$.
\end{enumerate}
\end{lemma}

\begin{proof}
We shall use the one-to-one correspondence from Definition \ref{defNC}, between cycles $C$
and harmonic functions $h_C$; and we first consider the cycle $C_0=\{0\}$, setting $h_0:= h_{C_0}$.
The fact that $h_0$ is continuous and harmonic follows from
\cite{DuJo05}. Also, from \cite[Theorem 4.1]{DuJo05} we know that
$$\chi_{\mathbf{N}_0}(\omega)=\lim_{k\rightarrow\infty}h_0(\tau_{\omega_k}\cdots \tau_{\omega_1}x),\quad\mbox{ for }P_x\mbox{ a.e.\ }\omega.$$
So, if $h_0$ has no zeroes then $\mathbf{N}_0$ must have $P_x$-measure $1$, otherwise there are some $\omega$ for which the limit is $0$, and, using the continuity of $h_0$, we obtain a contradiction.
The argument in the proof of the second assertion in the lemma is essentially the same, with the modification that we
now consider the union of the two subsets $\mathbf{N}_0$ and $\mathbf{N}_1$. Assuming that $h_0 + h_1$  has no zeroes,
then the same reasoning as above yields $\chi_{\mathbf{N}_0\cup\mathbf{N}_1}=1$ $P_x$ a.e., and therefore  $P_x(\mathbf{N}_0\cup\mathbf{N}_1)=1=h_0(x)+h_1(x)$.
\end{proof}

\begin{remark}  Our limit formula above for the indicator function
$v:= \chi_{\mathbf{N}_0}$ is a special case of a martingale limit argument
from \cite[Theorem 3.1]{DuJo05}, see also \cite{Gun00}. Our result
from [DuJo05] is further a close analogue for IFS-harmonic functions
of the familiar and classical Fatou-Primalov boundary-value theorem
in Fourier analysis \cite{Rud87}.

A special class of cocycles $v$ in \cite{DuJo05} are the indicator
functions of subsets $\mathbf{N}_C$ of the symbol space $\Omega$ which are
defined in terms of the shift on $\Omega$; i.e., $v = \chi_{\mathbf{N}_C}$.
The reader is referred to \cite{DuJo05} for additional details and to the list of cycles for the Bernoulli
convolutions above, i.e., the list before the statement of Theorem \ref{lemnocycles}.
\end{remark}

\par
We first need a definition and some technical results.
\begin{definition}\label{defcyclo}
Let $\lambda\in\br_+$, and suppose for some
$p\geq2$, $\lambda^p\in\mathbb{Q}$. We say that
$p$ is minimal if $\lambda^{p_1}\not\in\mathbb{Q}$ when $0<p_1<p$.
The number $p$ is called the minimal degree of $\lambda$.
\par
The minimal polynomial $m(x)\in\mathbb{Q}[x]$ of $\lambda$ is a
monic polynomial, $m(\lambda)=0$, of lowest degree with $\lambda$ as
root. It is irreducible in $\mathbb{Q}[x]$; see \cite[p.\
234]{Hun84}. Moreover, if $f(x)\in\mathbb{Q}[x]$ satisfies
$f(\lambda)=0$, then $m(x)$ divides $f(x)$.
\end{definition}
\begin{lemma}\label{lemcyclo0}
If $\lambda=q^{\frac1p}$ with $q\in\mathbb{Q}$ and $p\geq2$ minimal,
then the minimal polynomial of $\lambda$ is $x^p-q$.
\end{lemma}
\begin{proof}
Let $m(x)$ be the minimal polynomial of $\lambda$. Then $m(x)$ has
to divide $x^p-q$, therefore $m(x)$ has to be of the form
$$m(x)=\prod_{k}(x-\lambda z_k),$$
where the product is over some of the $p$'th roots of unity. Let
$1\leq l\leq p$ be the degree of $m$. The constant term is of the
form $\lambda^lz$ for some $z\in\bc$ with $|z|=1$. At the same time
it has to be rational. Thus $z=\pm 1$ and $\lambda^l$ is rational.
Since $p$ is minimal, it follows that $p=l$ and that $m(x)=x^p-q$.
\end{proof}
\begin{corollary}\label{lemcyclo}
Let $\lambda\in\br_+$, and let $p$ be its minimal degree, $p\geq2$.
Then there are no solutions
$$
\lambda^r=a\lambda^s+b,
$$
with $a,b\in\mathbb{Q}$, and $0\leq s<r<p$.
\end{corollary}
\begin{proof}
If $\lambda$ satisfies this condition, then the polynomial
$x^r-ax^s-b $ would be divisible by $x^p-q$; and this is impossible
since $r < p$.
\end{proof}

\begin{theorem}\label{proph-0=1}
Let $\lambda$ be in $(0,1)$.
\begin{enumerate}
\item
If $\lambda\not\in\{1-\frac1{2n}\,|\,n\in\bn\}$ then $h_0=1$.
\item
If $\lambda=1-\frac1{2n}$ then $h_0+h_1=1$.
\end{enumerate}
\end{theorem}

\begin{proof}
We use Lemma \ref{lemnozerois1}, and prove that $h_0$ has no zeroes. Using (\ref{eqh-0}), we see that it is enough to prove that, for each $x\in X_L$, there is a path $\omega$ that ends in $000\dots $ (or in $\frac14\frac14\frac14\dots $ for (ii)) which does not go through zeroes of $W$, i.e., $\tau_k\cdots \tau_1x$ is not a zero of $W$ for all $k\geq1$.
\par
Since $W=\cos^2\left(\frac{2\pi x}{\lambda}\right)$, the zeroes of $W$ are $Z:=\{\frac{\lambda}2(\frac12+n)\,|\,n\in\{0,1,\dots \}\}$.
Define now the set
$$A:=\bigcup_{p\geq 1}\lambda^{-p}Z.$$
If $x\not\in A$ then the desired path is $\omega=000\dots $.
\par
Assume now that $x\in A$. Take the smallest $p\geq 1$, such that
$\lambda^px\in Z$. There exists $m\in\{0,1.,\dots \}$ and $p\geq 1$
such that $x=\lambda^{-p}\frac{\lambda}{2}(\frac12+m)$. Since $p$ is
minimal, $\tau_0^kx\not\in Z$ for $k\leq p-1$. \par We consider
first the case when $\lambda\in(\smash{\frac12},1)$ is not of the form
$\lambda=\smash{\left(\frac{n+1/2}{m+1}\right)^{\frac1p}}$, for
$m,n\in\{0,1,\dots \}$, $p\geq1$. \par
We prove now that
$\tau_{1/4}\tau_0^{p-1}x\not\in A$. If not, then there exist $r\geq
1$ and $n\in\{0,1,\dots \}$ such that
$$\lambda\left(\lambda^{-1}\frac{\lambda}2\left(\frac12+m\right)^{\mathstrut}+\frac14\right)=\lambda^{-r}\frac\lambda2\left(\frac12+n\right).$$
This implies that $\lambda^r=\frac{n+1/2}{m+1}$, which contradicts the hypothesis.
 \par
Thus, in this case, the desired path is $\underbrace{00\dots 0}_{\!p-1\mbox{ times}\!}\frac1400\dots $.
\par
We can assume now that $\lambda$ is of the form
$\lambda=q^{\frac1p}$ for some $q$ rational and $p\geq 1$.
\par
The argument above shows that we may assume that $x\in
\lambda^{-1}Z$,
$x=\lambda^{-1}\frac{\lambda}2(n+\frac12)=\frac{n}2+\frac14$, for
some natural number $n\geq0$. Then $W(\tau_0x)=0$ so
$W(\tau_{1/4}x)=1$. We prove that $\tau_{1/4}\tau_{1/4}x\not\in
A\cup Z$.
\par
Suppose $\tau_{1/4}\tau_{1/4}x$ is in $A\cup Z$, then for some
$r\geq0$ and $m\in\{0,1,\dots \}$ we have
$$\lambda\left(\lambda\left(\left(\frac{n}2+\frac14\right)^{\mathstrut}+\frac14\right)^{\mathstrut}+\frac14\right)=\lambda^{-r}\frac{\lambda}2\left(m+\frac12\right),$$
so \begin{equation}\label{eqpoly2} \lambda^{r+1}\left(\frac
n2+\frac12\right)+\lambda^r\frac14=\left(\frac m2+\frac14\right). \end{equation}

We will show that this yields a contradiction with the fact that
$\lambda^p$ is rational for some $p\geq2$.
\par
We first need a definition and some technical results.

\par
Consider $p$ minimal such that $\lambda^p\in\mathbb{Q}$. If in
equation (\ref{eqpoly2}), $r+1\geq p$ then we can write $r+1=pi+j$
for $i\in\bn$ and $j<p$. Then
$\lambda^{r+1}=\lambda^{pi}\lambda^j\in\mathbb{Q}\lambda^j$,
therefore, using (\ref{eqpoly2}), we see that $\lambda$ satisfies an
equation of the form $\lambda^r=a\lambda^s+b$ with $0\leq s<r<p$,
and $a,b\in\mathbb{Q}$, and this contradicts Corollary
\ref{lemcyclo}.
\par
Finally we obtain that $\frac14\frac14000\dots $ is the desired path
that avoids zeroes of $W$.
\par
The remaining case is when $\lambda$ is rational. Take $x\in X_L$.
If $x=0$ then the path is $000\dots $. In (ii) if
$x=\frac{\lambda}{4(1-\lambda)}$, then the path is $\frac14\frac14\frac14\dots $.
\par
 Let $x\in
X_L\setminus\{0,\frac{\lambda}{4(1-\lambda)}\}$. We say that
$(\omega_1,\dots ,\omega_n)$ is a path that starts at $x$ and avoids
zeroes if
$W(\tau_{\omega_1}x),\dots ,W(\tau_{\omega_n}\cdots \tau_{\omega_1}x)$ are
all non-zero. Let $r_n$ be the number of paths that avoid zeroes.
Then $r_{n+1}\geq r_n$ because for any path $(\omega_1,\dots \omega_n)$
that avoids zeroes, since
$$W(\tau_0(\tau_{\omega_n}\cdots \tau_{\omega_1}x))+W(\tau_{1/4}(\tau_{\omega_n}\cdots \tau_{\omega_1}x))=1,$$
there is a $\omega_{n+1}\in\{0,1/4\}$ such that
$(\omega_1,\dots ,\omega_{n+1})$ is a path that avoids zeroes.
\par
Assume that $r_n$ is bounded. Then for some $n_0$, $r_{n+1}=r_n$ for all $n\geq n_0$.
Denote by
$$E_n:=\{\tau_{\omega_n}\cdots \tau_{\omega_1}x\,|\,(\omega_1,\dots ,\omega_n)\mbox{
is a path that avoids zeroes.}\}.$$ Since $r_n=r_{n+1}$ for $n\geq
n_0$, it means that for every $y\in E_n$, $n\geq n_0$, there is only
one choice to continue the path that starts at $x$ and ends at $y$
and avoids zeroes (otherwise $r_{n+1}\geq r_n+1$). Thus
$W(\tau_0y)=0$ or $W(\tau_{1/4}y)=0$.
\par
This implies that $y\in\tau_0^{-1}(Z)\cup \tau_{1/4}^{-1}(Z)$ where
$Z$ is the set of zeroes of $W$. Consequently $E:=\bigcup_{n\geq
n_0}E_n$ is a finite set.
\par
Now pick $y$ in $E$. We have that $\tau_0x$ or $\tau_{1/4}x$ is a
zero for $W$. Therefore  $W(\tau_{1/4}x)=1$ or $W(\tau_0x)=1$.
Let $\eta_1\in\{0,1/4\}$ such that $W(\tau_{\eta_1}y)=1$. Then
$\tau_{\eta_1}y\in E$. Inductively we can construct
$\eta_1,\dots ,\eta_m$ such that $W(\tau_{\eta_k}\cdots \tau_{\eta_1}x)=1$
and $\tau_{\eta_k}\cdots \tau_{\eta_1}x\in E$. Since $E$ is finite,
$\tau_{\eta_k}\cdots \tau_{\eta_1}x$ must end in a cycle. Moreover, the
previous argument shows that this will be a $W$-cycle.
\par
However the only cycles are the $1$-cycles $0$ and possibly
$\frac{\lambda}{4(1-\lambda)}$ in (ii). But note that for $z\in
X_L$, $\tau_{\alpha}z=0$ implies $z=0$, and
$\tau_{\alpha}z=\frac{\lambda}{4(1-\lambda)}$ implies
$z=\frac{\lambda}{4(1-\lambda)}$. Thus
$x\in\{0,\frac{\lambda}{4(1-\lambda)}\}.$
\par
Thus we have that $r_n$ is unbounded.
 Next we prove that, since $\lambda$ is rational, we cannot have
 $\tau_{\omega_n}\cdots \tau_{\omega_1}x=\tau_{\eta_n}\cdots \tau_{\eta_1}x$
 for $(\omega_1,\dots ,\omega_n)\neq(\eta_1,\dots ,\eta_n).$ Indeed, if
 they are equal then we obtain
 $$\omega_n+\lambda\omega_{n-1}+\dots +\lambda^{n-1}\omega_n=\eta_n+\lambda\eta_{n-1}+\dots +\lambda^{n-1}\eta_n.$$
 Then $\lambda$ will satisfy an equation of the form
 $$a_{r}\lambda^{r}+\dots +a_1\lambda+a_0=0,$$
 where $a_{i}\in\{0,1,-1\}$, not all of them $0$, and $r\geq 1$. Let $\lambda=a/b$ with $a,b$ mutually
 prime integers. We
 have$$a_{r}a^{r}+a_{r-1}a^{r-1}b+\dots +a_{0}b^r=0.$$
 But all the terms are divisible by $b$ except the first one, and
 this yields a contradiction.
 \par
 Since $r_n$ is unbounded and, since for different paths of the same
 length starting at $x$ the endpoints are different, it follows that
 the cardinality of the set $E_n$ increases to infinity.
 \par
 However, note that the set $A$ is finite because $Z$ is finite and
 $\lambda^{-r}z\not\in X_L$ for $r$ large. It follows then, that if
 we choose $n$ such that the cardinality of $E_n$ is bigger then the
 one of $A$, there is a path of length $n$ that starts at $x$, avoids zeroes, and
 ends outside $A$. Hence this path continued with $000\dots $ is the
 desired one.
\end{proof}
\begin{corollary}\label{cor2-1}
Consider the Bernoulli IFSs for all $\lambda$ in the open interval $(0, 1)$. For $\lambda$ fixed,
and for all $x$, each of the measures $P_x$ on $\Omega$ is countably atomic. By this we mean
that each $P_x$ is supported on a certain countable subset $S$ of $\Omega$, and each point $\omega$ in $S$ is an atom for
$P_x$.
\end{corollary}

\begin{remark}
This property from the conclusion of the corollary for the process $P_x$ is called ``tightness''. See, e.g., \cite{Gun00}, \cite{Jor06}, and \cite{Fal93},
where an analogous process $(\Omega, P_x)$ is studied in a different context: orthogonality relations for the wavelet scaling function.
Our result, Corollary \ref{cor2-1} for the Bernoulli convolutions, to the effect that the measures
$P_x$ are purely atomic, stands in contrast with a different result below, for the
Bernoulli measures $\nu_\lambda$. In fact, we show in Corollary \ref{cornoatombern} that the strongly invariant
measures for affine IFSs (which includes the Bernoulli measures $\nu_\lambda$) are all non-atomic.
\end{remark}

Let $\lambda\in(0,1)$ and consider the random variable $\sum_{j=0}^\infty(\pm1)\lambda^j$, i.e., independently
distributed ones and minus ones. It follows from the previous sections that the distribution of this random variable
is the measure $\nu_\lambda=p_{(\frac12,\frac12)}\circ\pi_\lambda^{-1}$.

\begin{lemma}\label{lemnuhat}
Let $\lambda\in(0,1)$, and let $\nu_\lambda$ be the corresponding infinite convolution IFS measure. Then the Fourier transform
\begin{equation}\label{eqnuhat1}
\hat\nu_\lambda(t):=\int e^{2\pi ixt}\,d\nu_\lambda(x)
\end{equation}
satisfies
\begin{equation}\label{eqnuhat2}
\hat\nu_\lambda(t)=\prod_{n=0}^\infty\cos(2\pi\lambda^nt).
\end{equation}
\end{lemma}

\begin{proof}
By Lemma \ref{lembernmeas}, (\ref{eqnu3-3}), we get
\begin{equation}\label{eqnuhat3}
\nu_\lambda=\frac12\left(\nu_\lambda\circ\tau_+^{-1}+\nu_\lambda\circ\tau_-^{-1}\right),
\end{equation}
or equivalently,
$$\int f(x)\,d\nu_\lambda(x)=\frac12\left(\int f(\tau_+(x))\,d\nu_\lambda(x)+\int f(\tau_-(x))\,d\nu_\lambda(x)\right),$$
see (\ref{eqtaupm}).
Using this and (\ref{eqnuhat1}), we get
$$\hat\nu_\lambda(t)=\int e^{2\pi ixt}\,d\nu_\lambda(x)$$$$=\frac12\left(\int e^{2\pi i(\lambda x+1)t}\,d\nu_\lambda(x)+\int e^{2\pi i(\lambda x-1)t}\, d\nu_\lambda(x)\right)$$$$=\frac{e^{2\pi it}+e^{-2\pi it}}{2}\int e^{2\pi i\lambda xt}\,d\nu_\lambda(x)=\cos(2\pi t)\hat\nu_\lambda(\lambda t)$$$$=\dots 
=\prod_{k=0}^{n-1}\cos(2\pi\lambda^kt)\hat\nu_\lambda(\lambda^nt),$$
by iteration.
\par
A limit argument, using that $$\lim_{n\rightarrow\infty}\hat\nu_\lambda(\lambda^nt)=\hat\nu_\lambda(0)=1,$$
now yields the desired infinite-product formula (\ref{eqnuhat2}).
\end{proof}

We will use the following notation: for a finite word $\omega=\omega_1\dots \omega_n\in\{0,1/4\}^n$, we denote by $W^{\omega}$ and $\tau_\omega$ the following objects:
$$W^{\omega}(x):=W(\tau_{\omega_n}\cdots \tau_{\omega_1}x)\cdots W(\tau_{\omega_1}x),\quad W^\emptyset(x)=1,$$
$$\tau_\omega x=\tau_{\omega_n}\cdots \tau_{\omega_1}x,\quad \tau_{\emptyset}x=x.$$

\begin{theorem}\label{thunique}
Let $\lambda\in(0,1)$, $\lambda$ not of the form $\lambda=1-\frac1{2n}$, $n\in\bn$.
Let
$$\mathcal{W}_0:=\{\emptyset\}\cup\{\omega=\omega_1\dots \omega_n\,|\,\omega_i\in\{0,1/4\}, \omega_n=1/4, n\geq1\}.$$
The function $|\hat\nu_\lambda|^2$ is the only continuous function of $X_\lambda$ that satisfies the following identity
\begin{equation}\label{eqmainid}
\sum_{\omega\in\mathcal{W}_0}W^\omega(x)f(\tau_\omega x)=1,\quad(x\in X_L).
\end{equation}
\end{theorem}

\begin{proof}
First we prove that $|\nu_\lambda|^2$ satisfies (\ref{eqmainid}). From Theorem \ref{proph-0=1} we know that
\begin{equation}\label{eqmainid1}
P_x(\mathbf{N}_0)=1.
\end{equation}
Take $\omega\in\mathbf{N}_0$. If $\omega=000\dots $, then, from \cite{DuJo05}, we know that
$$P_x(\{000\dots \})=\prod_{k=1}^\infty W(\tau_0\cdots \tau_0x)=|\hat\nu_\lambda(x)|^2,$$
where we used Lemma \ref{lemnuhat} in the last equality.
\par
If $\omega\neq000\dots $, then $\omega=\omega_1\dots \omega_n00\dots $, where $\omega_n=1/4$. Then
$$P_x(\{\omega\})=W(\tau_{\omega_1}x)\cdots W(\tau_{\omega_n}\cdots \tau_{\omega_1}x)\prod_{k=1}^\infty W(\tau_0^k\tau_{\omega_n}\cdots \tau_{\omega_1}x)$$
$$=W^{\omega_1\dots \omega_n}(x)|\hat\nu_\lambda(\tau_\omega x)|^2.$$
Thus equation (\ref{eqmainid1}) can be rewritten as (\ref{eqmainid}).
\par
Take now a continuous function $f$ on $X_L=[0,\frac{\lambda}{4(1-\lambda)}]$ which satisfies (\ref{eqmainid}). We rewrite this identity, by enumerating the finite words that end in $1$:
\begin{multline*}
1=f(x)+W(\tau_{1/4}x)f(\tau_{1/4}x)+\cdots \\
{}+\sum_{\omega_1,\dots ,\omega_{n-1}}W(\tau_{\omega_0}x)\cdots W(\tau_{\omega_{n-1}}\cdots \tau_{\omega_1}x)\qquad\qquad\\
{}\cdot W(\tau_{1/4}\tau_{\omega_{n-1}}\cdots \tau_{\omega_1}x)f(\tau_{1/4}\tau_{\omega_{n-1}}\cdots \tau_{\omega_1}x)+\cdots.
\end{multline*}
Therefore
\begin{equation}\label{eqmainid2}
1=f(x)+R_W^0(W\circ\tau_{1/4}\,f\circ\tau_{1/4})(x)+\dots +R_W^{n-1}(W\circ\tau_{1/4}\,f\circ\tau_{1/4})(x)+\cdots.
\end{equation}
In particular, $R_W^{n}(W\circ\tau_{1/4}\,f\circ\tau_{1/4})$ converges pointwise to $0$.
\par
We apply the operator $I-R_W$ to (\ref{eqmainid2}). (More precisely, we consider equation (\ref{eqmainid2}) at $\tau_0x$ and $\tau_{1/4}x$, and multiply by
$W(\tau_0x)$ and $W(\tau_{1/4}x)$ respectively. Then we subtract both these equations from (\ref{eqmainid2}). All the series involved in the computation converge pointwise). We obtain
$$0=f-R_Wf+W\circ\tau_{1/4}\,f\circ\tau_{1/4}-R_W(W\circ\tau_{1/4}\,f\circ\tau_{1/4})+\cdots $$$${}+R_W^{n-1}(W\circ\tau_{1/4}\,f\circ\tau_{1/4})-R_W^n(W\circ\tau_{1/4}\,f\circ\tau_{1/4})+\cdots $$
$$=f-R_Wf+W\circ\tau_{1/4}\,f\circ\tau_{1/4}=f-W\circ\tau_0 f\circ\tau_0.$$
Thus $f$ satisfies the equation
\begin{equation}\label{eqmainid3}
f(x)=W(\lambda x)f(\lambda x),\quad x\in\left[0,\frac{\lambda}{4(1-\lambda)}\right].
\end{equation}
Note also that equation (\ref{eqmainid2}) implies that $f(0)=1$. Then, using the relations in the proof of Lemma \ref{lemnuhat}, we see that
the function $g:=f/|\hat\nu_\lambda|^2$ verifies the equation
$$g(x)=g(\lambda x),$$
for all $x$ in $[0,\frac{\lambda}{4(1-\lambda)}]$ except when $x$ is a zero for $\hat\nu_\lambda$. The zeroes of $\hat\nu_\lambda$ inside the interval $X_L$ are finite in number, because they are of the
form $\lambda^{-k}z$ where $z$ is a zero of $W$.
\par
In addition, the function $g$ is continuous outside this finite set; it is also continuous at $0$ and $g(0)=1$. This implies
that $g(x)=g(\lambda^nx)=\lim_{n\rightarrow\infty}g(\lambda^nx)=g(0)=1$ for $x$ outside some finite set. Finally we can conclude that $g=1$ except at some finite number of points, and this implies that
$f=|\hat\nu_\lambda|^2$.
\end{proof}

\begin{proposition}\label{propunifconv}
The convergence in \textup{(\ref{eqmainid})} is uniform on $X_\lambda=[0,\frac{\lambda}{4(1-\lambda)}]$.
\end{proposition}

\begin{proof}
By Theorem \ref{thunique}, we have that $f=|\hat\nu_\lambda|^2$. But this is a non-negative continuous function and the same is true for $W$. Thus the uniform convergence follows from the pointwise convergence and Dini's theorem.
\end{proof}

\section{$\lambda$ representations}\label{reps}

While there is already a substantial literature on the structure of infinite Bernoulli convolutions,
we take here a different approach. In this paper, the systems $(X_\lambda, \nu_\lambda)$  serve as a link between
the class of iterated function systems (IFS) and the more general scale-similarity notions in analysis, e.g.,
the scale covariance which is intrinsic to wavelet theory; see the discussion around equations (\ref{eqinsdz1}) and (\ref{eqinsdz2}) above.
We showed in \cite{DuJo05} that in general IFSs carry a discrete harmonic analysis which is based on $W$-cycles.
For general IFSs, the family of $W$-cycles may in fact be quite complicated and not especially computable.
Our present Theorem \ref{lemnocycles} in the previous section shows that for each of the Bernoulli systems $X_\lambda$,
we can have at most two $W$-cycles. We found them, and worked out the harmonic analysis. This harmonic analysis as
a result is somewhat coarse.

          The purpose of this section is to relate this theorem and its implications to the geometry
          of the Bernoulli systems $X_\lambda$ and to their encoding mappings  $\pi_\lambda \colon  \Omega\rightarrow X_\lambda$.
          While it is possible that some of our results spelled out in Proposition \ref{lembern1} below are known,
          we have not been able to locate them in the literature.

 For a parameter $\lambda$ in the open interval $(0, 1)$,  consider the random
variable $\sum_k(\pm\lambda^k)$, with independent $\frac12$-$\frac12$-distribution of the $\pm$ signs in the series, and let $\nu_\lambda$ be the corresponding distribution measure. Solomyak \cite{So95} showed that
for Lebesgue a.e.\ $\lambda$ in $(1/2, 1)$,  $\nu_\lambda$ has an $L^2$ density. Nonetheless, when $\lambda > 1/2$,  preciously
little is known about the explicit properties of this one-parameter family of measures $\nu_\lambda$. A fundamental result is due to Paul Erd\"os \cite{Er39} who
showed that $\nu_\lambda$ is purely singular when the reciprocal $\lambda^{-1}$ is a Pisot number, i.e., an algebraic number whose Galois conjugates are strictly less than one in absolute value.
While P. Erd\"{o}s' theorem \cite{Er39} stands out, it has now been extended to a few larger classes of algebraic numbers
(see especially the deep results by Garcia \cite{Ga62} about absolute continuity, and Kahane \cite{Kah71}
proving that the Fourier transform of $\nu_\lambda$ does not vanish at infinity when $\lambda^{-1}$ is a Salem number), yet
 there is little known about specific values of $\lambda$; for example, few explicit properties are known about
 the measure $\nu_{2/3}$. Yet, the isolated results in the literature do depend on certain algebraic/arithmetic
 properties of the number $\lambda$. In fact the Diophantine properties of $\lambda$ will play a role in our present analysis.


In order to apply our Fourier duality to the classical $\lambda$-systems, we need a few geometric facts about these systems. We have therefore assembled them in the next proposition. While they may be known, we have not been able to find them exactly in the form that is needed for our present harmonic analysis. So we include them without proof for the benefit of the reader.

\begin{proposition}\label{lembern1}
Let $\lambda\in(0,1)$ and $B=\{0,1\}$, and let $(\tau_b)$, $X_B$ be the corresponding affine IFS.
\begin{enumerate}
\item[(a)] If $\lambda<1/2$, then $X_B$ is a fractal.
\item[(b)] If $\lambda\geq1/2$, then
\begin{equation}\label{eqlembern16}
X_B=\left[0,\frac{\lambda}{1-\lambda}\right],
\end{equation}
i.e., a closed interval.
\item[(c)] If $0<\lambda<1/2$, then $\pi_\lambda\colon \Omega\rightarrow X_B$ is a bijection.
\item[(d)] If $\frac{1}{2}<\lambda<1$, then for every point $x\in (0,\frac{\lambda}{1-\lambda})$, $\pi_\lambda^{-1}(x)$ is uncountable.
\end{enumerate}
\end{proposition}

\begin{proof}
See \cite{JoPe96}, \cite{Edg98}, \cite{JeWi35}, \cite{KeWi35}, \cite{AnIl01} and the references there\-in. 
\end{proof}
\begin{corollary}\label{corhausdim}
If $0<\lambda<\frac12$, then $X_B$ is a Cantor set of Lebesgue measure zero and Hausdorff dimension $d_H=\frac{\ln(2)}{\ln(\lambda^{-1})}$.
\end{corollary}
\begin{proof}
We first show that the Lebesgue measure of the support set $X_\lambda$ is
a fractal, and that its Lebesgue measure is zero.
Start with the interval $[0,\frac{\lambda}{1-\lambda}]$ and omit the middle open interval $(\frac{\lambda^2}{1-\lambda},\lambda)$ of length $\frac{\lambda(1-2\lambda)}{1-\lambda}$. Continue to omit the middle open intervals at each $\lambda$-scale.
\par
The sum of all the omitted middle fractions is
$$\frac{\lambda(1-2\lambda)}{1-\lambda}(1+2\lambda+(2\lambda)^2+\cdots )=\frac{\lambda(1-2\lambda)}{1-\lambda}\frac{1}{1-2\lambda}=\frac{\lambda}{1-\lambda},$$
so the omitted intervals have full measure in the interval $[0,\frac{\lambda}{1-\lambda}]$. The stated formula for the Hausdorff dimension $d_H$ now follows from
the discussion above, and from standard scaling theory; see, e.g.,
\cite{Edg98}.
\end{proof}

\begin{remark}\label{remfax1}
For $y\in\br$, set $\langle y\rangle :=y-\lfloor y\rfloor=$the fractional part of $y$. Then $0\leq\langle y\rangle <1$. When $0<\lambda<1$ is given, set $r_\lambda\colon [0,1)\rightarrow[0,1)$,
\begin{equation}\label{eqrlambda}
r_\lambda(y):=\langle \lambda^{-1}y\rangle 
\end{equation}
We note that this is a dynamical system. This system has uses in symbolic dynamics, for example in the study of subshifts, see \cite{Wal78}.

We now recall a recursive number-representation based on (\ref{eqrlambda}). It is simply a recursive algorithm for a
$\lambda$-representation of certain real numbers $x$. When $\lambda<1/2$, it
does yield the representation of points in the fractal $X_\lambda$, but
if $\lambda\geq 1/2$  (Lemma \ref{lemfax2}), this representation typically has two
defects: first, the required alphabet $A$ is then bigger than $\{0,1\}$;
and secondly, the dynamical systems expansion (\ref{eqomr}) only
offers {\it one} representation of the admissible numbers $x$. Our result
Proposition \ref{lembern1} above in fact keeps track of the number of different
representations for the same $x$; and the alphabet $A$ in Proposition \ref{lembern1}
is fixed, $A =\{0,1\}$.
\par
When $0<\lambda<\frac12$, we saw that the expansion $x=\sum_{k=1}^\infty\omega_k\lambda^k$, is unique, $\omega_k\in\{0,1\}$ for
all $x\in X_B$. An application of (\ref{eqrlambda}) shows that
\begin{equation}\label{eqomr}
\omega_k=\left\lfloor\frac{1}{\lambda}r_\lambda^{k-1}(x)\right\rfloor.
\end{equation}
If $0<\lambda<\frac12$, then there are only two possibilities, $0$ or $1$, for $\omega_k$. Not so if $\frac12\leq\lambda<1$.
\par
There is a more general expansion $x=\sum_{i=1}^\infty\omega_i\lambda^i$ which has the $\pi_\lambda$ expansion as a special case when $\lambda<\frac12$.
\par
For $y\in\br$, let $\lfloor y\rfloor=$integral part, $\langle y\rangle :=y-\lfloor y\rfloor=$fractional part.
\par
If $0<\lambda<1$, set $a=$ the largest integer with $a<\frac{1}{1-\lambda}$. Set
$$r_\lambda(y):=\left\langle\frac{1}{\lambda}y\right\rangle ,\quad\omega_k:=\left\lfloor\frac{1}{\lambda}r_\lambda^{k-1}(x)\right\rfloor\in\{0,1,\dots ,a\}.$$

\end{remark}

\begin{lemma}\label{lemfax2}
Let $0<\lambda<1$, and let $A$ be the set of integers $\{0,1,\dots ,a\}$ where $a<\frac{1}{1-\lambda}$. (Note that if $\lambda<\frac12$, then $A=\{0,1\}$.)
Then every $x\in X_B$ has a convergent expansion
\begin{equation}\label{eqx=omega-k}
x=\sum_{k=1}^\infty\omega_k\lambda^k,
\end{equation}
with $\omega_k\in A$ as in (\ref{eqomr}).
\end{lemma}
\begin{proof}
We only need to show that (\ref{eqx=omega-k}) is convergent. But the
recursive algorithm shows that
$$0\leq x-\sum_{j=1}^k\omega_j\lambda^j\leq a\frac{\lambda^k}{1-\lambda},$$
where $a$ is the largest integer $<\frac{1}{1-\lambda}$.
Since $\lambda<1$, the result follows.

\end{proof}

The following simple lemma shows that the IFS associated to $\lambda\in(0,1)$ and an arbitrary pair $B=\{a,b\}$ is equivalent up to some scaling and translation to
the one associated to $\{0,1\}$. Therefore in analyzing these systems, we can take any pair $B$.
\begin{lemma}
Let $\lambda\in(0,1)$ and $B=\{a,b\}$ with $a<b$ in $\br$. Then the affine function $\alpha(x)=(b-a)x+a\frac{\lambda}{1-\lambda}$ has the property that, for $(a_i)_i\in\{0,1\}^\bn$,
$$\alpha(\sum_{i=1}^\infty a_i\lambda^i)=\sum_{i=1}^\infty b_i\lambda^i,\mbox{ where }b_i:=\left\{\begin{array}{ccc}
a&\mbox{ if }&a_i=0,\\
b&\mbox{ if }&a_i=1.\end{array}\right.$$
\end{lemma}
\begin{proof}
It is enough to see that $$a\sum_{i=1}^\infty a_i\lambda^i+(b-a)\frac\lambda{1-\lambda}=\sum_{i=1}^\infty (aa_i+(b-a))\lambda^i.$$
\end{proof}

\section{Strongly invariant measures}\label{meas}
Let $(X,(\tau_i)_{i=1}^N)$ be an IFS, and suppose it has an encoding mapping $\pi\colon \Omega\rightarrow X$ which is measurable with respect to two $\sigma$-algebras
\begin{itemize}
\item $\mathfrak{B}(\Omega)$: the $\sigma$-algebra on $\Omega$ which is generated by the cylinder sets.
\item $\mathfrak{B}(X)$: the Borel $\sigma$-algebra on $X$.
\end{itemize}
For every $(p_i)_{i=1}^N$, $p_i>0$, $\sum_{i=1}^Np_i=1$, there is a unique probability measure $\mu_{(p)}$, on $\Omega$, called the Bernoulli measure, and determined by
\begin{equation}\label{eqbernmeas1}
\mu_{(p)}(\omega_1\omega_2\dots \omega_n\times\bz_N\times\bz_N\times\cdots )=p_{\omega_1}p_{\omega_2}\cdots p_{\omega_n}.
\end{equation}

\begin{lemma}\label{lembernmeas} Let $(X,(\tau_i))$ be an IFS which has a measurable encoding mapping $\pi$.
Let $(p_i)$ be as above, and let $\mu_{(p)}$ be the Bernoulli measure on $\Omega$. Then the measure
\begin{equation}\label{eqnu3-2}
\nu_{(p)}^X:=\mu_{(p)}\circ\pi^{-1}
\end{equation}
on $X$ satisfies
\begin{equation}\label{eqnu3-3}
\nu_{(p)}^X=\sum_{i=1}^Np_i\nu_{(p)}^X\circ\tau_i^{-1},
\end{equation}
and
$$\mbox{supp}(\nu_{(p)}^X)=X,$$
where $\mbox{supp}$ denotes the support.
\end{lemma}

\begin{proof}
The measure in (\ref{eqnu3-2}) is defined by $\mu_{(p)}(\pi^{-1}(E))$ for $E\in\mathfrak{B}(X)$, where
$$\pi^{-1}(E)=\{\omega\in\Omega\,|\,\pi(\omega)\in E\};$$
and similarly
$$(\nu_{(p)}^X\circ\tau_i^{-1})(E)=\nu_{(p)}^X(\tau_i^{-1}(E)).$$
It follows from the definition of $\mu_{(p)}$ in (\ref{eqbernmeas1}) that
$$\mu_{(p)}=\sum_{i=1}^Np_i\mu_{(p)}\circ\sigma_i^{-1},$$
where $\sigma_i\colon \Omega\rightarrow\Omega$ is defined by $\sigma_i(\omega_1\omega_2\dots ):=(i\omega_1\omega_2\dots )$. From (\ref{eqdef33}), we conclude that
$$\pi\circ\sigma_i=\tau_i\circ\pi.$$
\par
Let $E\in\mathfrak{B}(X)$. Then
$$\nu_{(p)}^X(E)=\mu_{(p)}(\pi^{-1}(E))=\sum_{i=1}^Np_i\mu_{(p)}(\sigma_i^{-1}\pi^{-1}(E))$$
$$=\sum_{i=1}^Np_i\mu_{(p)}(\pi^{-1}\tau_i^{-1}(E))=\sum_{i=1}^Np_i\nu_{(p)}^X(\tau_i^{-1}(E)),$$
which is the desired formula (\ref{eqnu3-3}).
\par
We now prove that $\nu_{(p)}^X$ has full support. The proof is indirect: Suppose $\emptyset\neq G\subset X$ is open and $\nu_{(p)}^X(G)=0$.
Then there is a finite word $(i_1\dots i_n)$ such that
$$\sigma_{i_1}\cdots \sigma_{i_n}(\Omega)\subset\pi^{-1}(G).$$
But then
$$\nu_{(p)}^X(G)=\mu_{(p)}(\pi^{-1}(G))\geq\mu_{(p)}(\sigma_{i_1}\cdots \sigma_{i_n}\Omega)=p_{i_1}\cdots p_{i_n},$$
contradicting that $p_i>0$ for all $i$. Hence $\nu_{(p)}^X$ does not vanish on non-empty open subsets in $X$, which is to say
$$\mbox{supp}(\nu_{(p)}^X)=X;$$
see \cite[page 58]{Rud87}.
\end{proof}
\begin{definition}\label{defstronglyinv}
A probability measure $\nu_{(p)}^X$ on $X$ satisfying equation (\ref{eqnu3-3}) is called strongly invariant.
\end{definition}
Lemma \ref{lembernmeas} proves in particular that strongly invariant measures exist. When the maps $\tau_i$ are contractions
the strongly invariant measure is also unique (see \cite{Hut81}).
\par
We show that a class of these strongly invariant measures are non-atomic, i.e., that
$$\nu_{(p)}^X(\{x\})=0,\quad\mbox{for }x\in X.$$
\begin{definition}
Let $(X, (\tau_i)_{i=1,\dots ,N})$ be an IFS, and let $x$ be in $X$. Then the backward orbit of $x$, $O^-(x)$ is defined as
$$O^-(x):=\bigcup\{\tau_{\omega_n}^{-1}\cdots \tau_{\omega_1}^{-1}x\,|\,n\in\bn, \omega_1,\dots ,\omega_n\in\{1,\dots ,N\}\}.$$
\end{definition}
\par
Our result in this section is based on an {\it a priori} observation about atoms for strongly invariant measures $\nu$
(unique or not!) for IFSs. In this generality, our result states that if an IFS with one-to-one maps has atoms, it has a maximal atom,
 i.e., an atom $x$ such that  $a:=\nu(\{x\})\geq\nu(\{y\})$ for all other atoms $y$.  We then show that all the non-empty
 sets in the backward orbit $O^-(x)$ from a maximal atom $x$ consist of atoms of the same measure. We note that this fact
 strongly restricts those IFS that can have atoms; and we prove in Corollary \ref{cornoatombern} that all contractive
 affine IFS are non-atomic.

\begin{theorem}\label{thnoatoms}
Let $(X,(\tau_i)_{i\in\{1,\dots ,N\}})$ be an iterated function system, $(p_i)_{i=1}^N$ given as above, with $p_i>0$ for all $i$, and let $\nu:=\nu_{(p)}$
be the unique strongly invariant probability measure. Assume that for all $x\in X$, either the backward orbit of $x$
is infinite, or there are some $\omega_1,\dots ,\omega_n$ such that
$\tau_{\omega_n}^{-1}\cdots \tau_{\omega_1}^{-1}(x)=\emptyset$.
Then the measure $\nu$ has no atoms.
\end{theorem}

\begin{proof}
Since $\nu$ is strongly invariant, it satisfies the equation
\begin{equation}\label{eqatom}
\nu(\{x\})=\sum_{k=1}^Np_k\nu(\tau_k^{-1}x)\quad(x\in X).
\end{equation}
(Since $\tau_i$ are contractions hence injective, each set $\tau_k^{-1}x$ contains at most one point.)
Since $\nu$ is a finite measure, there are only countably many atoms and they can be listed in a sequence that converges to zero. Therefore there is a biggest atom, i.e., one with maximal measure. Let $x_0$ be an atom such that $\nu(\{x_0\})\geq\nu(\{x\})$ for all $x\in X$.
\par
The sets $\tau_k^{-1}x_0$ have measure less than $\nu(\{x_0\})$ by maximality.
From equation (\ref{eqatom}) we see that none of the sets $\tau_k^{-1}x_0$ can have measure strictly less
than $\nu(\{x_0\})$. Thus $\nu(\tau_k^{-1}x_0)=\nu(\{x_0\})$, for $k\in\{1,\dots ,N\}$.
By induction, we obtain that $\nu(\tau_{\omega_n}^{-1}\cdots \tau_{\omega_1}^{-1}x_0)=\nu(\{x_0\}),$
for all $n\in\bn$ and all $\omega_1,\dots ,\omega_n\in\{1,\dots ,N\}$. In particular, all sets $\tau_{\omega_n}^{-1}\cdots \tau_{\omega_1}^{-1}x$ are non-empty. But, using the hypothesis,
we obtain that there is an infinite number of points that have measure equal to $\nu(\{x_0\})>0$, and
this contradicts the fact that $\nu$ is finite.

\end{proof}

\begin{lemma}\label{leminfinvorb}
Let $(\tau_b)_{b\in B}$ be an affine IFS on $\br^d$ as before with $\#B\geq2$.
If $x\in\br^d$ then, for $n\in\bn$ the set $\{\tau_{\omega_n}^{-1}\cdots \tau_{\omega_1}^{-1}x\,|\,\omega_1,\dots ,\omega_n\in B\}$ has at least $n$ elements.
\par
Also the set
$\{\tau_{\omega_n}\cdots \tau_{\omega_1}x\,|\,\omega_1,\dots ,\omega_n\in B\}$ has at least $n$ elements.
\end{lemma}
\begin{proof}
First, note that $\tau_{\omega_n}^{-1}\cdots \tau_{\omega_1}^{-1}x=\tau_{\eta_n}^{-1}\cdots \tau_{\eta_1}^{-1}x$ if and only if
$$\omega_1+\dots +R^{n-1}\omega_n=\eta_1+\dots +R^{n-1}\eta_n.$$
Then, take $b_1\neq b_2$ in $B$, and for $i\in\{1,\dots ,n\}$, define $$\omega_k^{(i)}:=\left\{\begin{array}{ccc}
b_1,&\mbox{if}&k=i,\\
b_2,&\mbox{if}&k\neq i.\end{array}
\right.$$
Then we cannot have for $i>j$, $$\sum_{k=1}^nR^{k-1}\omega_k^{(i)}=\sum_{k=1}^nR^{k-1}\omega_k^{(j)},$$
because that implies $(R^{i-1}-R^{j-1})(b_1-b_2)=0$ so $R^{i-j}(b_1-b_2)=(b_1-b_2)$, ($R^{j-1}$ is invertible).
But this is false since $R$ is expansive so $R^{i-j}$ does not have the eigenvalue $1$.
We used that the scaling matrix $R$ is invertible. But $b_1-b_2$ cannot be an eigenvector of $R^{i-j}$
since $1$ is not in the spectrum of $R$.
Thus the points $\tau_{\omega_n^{(i)}}^{-1}\cdots \tau_{\omega_1^{(i)}}^{-1}x$ are distinct,
and this proves the first assertion. The same argument works for the second statement.
\end{proof}

\begin{corollary}\label{cornoatombern}
The strongly invariant measure of an affine IFS has no atoms.
\end{corollary}
\begin{proof}
Follows directly from Theorem \ref{thnoatoms} and Lemma \ref{leminfinvorb}.
\end{proof}

Next, we continue our study from \cite{DuJo03} of infinite
convolutions, and the corresponding infinite product formulas for
the associated Hutchinson measures.

The Fourier transform of $\nu$ may be used for discerning properties of $\nu$ itself, for example via Wiener's test on the
mean-square Cesaro sum. We will illustrate this in the next remark for the Bernoulli measures $\nu_\lambda$. We show
that the mean-square Cesaro sum is zero of high order. So this conclusion is stronger than $\nu_\lambda$ being non-atomic.

\begin{remark}
We can use a criterion of N. Wiener for the existence of atoms: see, e.g., \cite{Kat04}. Wiener's theorem states that
the following Cesaro limit exists and yields the sum absolute square of the atoms, i.e.,
\begin{equation}\label{eqwie1}
\lim_{T\rightarrow\infty}\frac{1}{2T}\int_{-T}^T|\hat\nu(t)|^2\,dt=\sum_{x\in\{\mbox{atoms}\}}|\nu(\{x\})|^2.
\end{equation}
In particular, if the expression on the left-hand side in (\ref{eqwie1}) is zero, then there are no atoms. For the Bernoulli systems, we will show that this limit on
the left is zero of exponential order.
\par
Set $$s(L)=\frac{1}{2L}\int_{-L}^L|\hat\nu(t)|^2\,dt,$$ for $L\in\br_+$. An application of the formula for $\hat\nu$ in the
proof of Lemma \ref{lemnuhat} now yields
\begin{equation}\label{eqwie2}
s(\lambda^{-n}T)=\frac{1}{2T}\int_{-T}^T\prod_{k=1}^n\cos^2(2\pi\lambda^{-k}t)|\hat\nu(t)|^2\,dt,
\end{equation}
for all $T\in\br_+$, and all $n\in\bn$. In particular, the sequence is monotone
$$s(\lambda^{-n}T)\leq\dots \leq s(\lambda^{-1}T)\leq s(T).$$
Therefore the limit $\lim_{n}s(\lambda^{-n}T)$ exists. \par
Note that, since $\lambda\in(0,1)$, $\lambda^{-n}T\rightarrow\infty$ as
$n\rightarrow\infty$.
So fix $T$, and let $n$ tend to infinity in (\ref{eqwie2}): This limit therefore coincides with the limit (\ref{eqwie1}). In fact, we show below that the limit in (\ref{eqwie2}) is $\mathcal{O}(2^{-n})$.
\par
Recall further from Example \ref{exbern2}, that $\nu=\nu_\lambda$ is supported in the closed interval $[-\frac1{1-\lambda},\frac1{1-\lambda}]$,
so the transform $t\mapsto\hat\nu(t)$ is analytic, in particular $C^\infty$.
\par
The product under the integral sign in (\ref{eqwie2}) is
\begin{equation}\label{eqwie3}
\prod_{k=1}^n\cos^2(2\pi\lambda^{-k}t)=2^{-n}\prod_{k=1}^n(1+\cos(4\pi\lambda^{-k}t)).
\end{equation}
But $\prod_{k=1}^n(1+\cos(4\pi\lambda^{-k}t))$ is a Riesz-product, which by \cite[Sections 1.3-1.4]{Kat04} represents a
distribution $D$ of compact support. In particular
$$\lim_{n\rightarrow\infty}\int_{-T}^T\prod_{k=1}^n(1+\cos(4\pi\lambda^{-k}t))|\hat\nu(t)|^2\,dt$$
is finite, i.e., it is $D(|\hat\nu|^2)$.
\par
Returning to (\ref{eqwie3}), we conclude that
$$\lim_{n\rightarrow\infty}s(\lambda^{-n}T)=\mathcal{O}(2^{-n})=0.$$
Substituting back into Wiener's formula (\ref{eqwie1}), the result follows.
\end{remark}

\begin{acknowledgements}The co-authors are pleased to acknowledge helpful
discussions with Professors  Keri Kornelson, David Larson, K.-S. Lau, Kathy
Merrill, Paul Muhly, Roger Nussbaum, Judy Packer, Steen Pedersen, and Yang Wang.
\end{acknowledgements}


\begin{thebibliography}{ccccccc}


\bibitem[ALTW04]{ALTW04}  Aldroubi, Akram; Larson, David; Tang, Wai-Shing; Weber, Eric Geometric
aspects of frame representations of abelian groups. Trans. Amer. Math. Soc. 356
(2004), no. 12, 4767--4786.

\bibitem[AnIl01]{AnIl01}
R. Anisca, M. Ilie, {\em A technique of studying sums of central Cantor sets}.
Canad. Math. Bull.  44  (2001),  no. 1, 12--18.

\bibitem[Ba00]{Ba00}
 V. Baladi, {\em Positive transfer operators and decay of correlations},
World Sci, Singapore 2000.
\bibitem[BKS91]{BKS91}
T. Bedford, M. Keane and C. Series, eds.; {\em Ergodic theory, symbolic dynamics, and hyperbolic spaces. Papers from the Workshop on Hyperbolic Geometry and Ergodic Theory held in Trieste, April 17--28, 1989}. Oxford Science Publications. The Clarendon Press, Oxford University Press, New York, 1991. xvi+369 pp.
\bibitem[BoPa05]{BoPa05}
Bodmann, Bernhard G.; Paulsen, Vern I. {\em Frames, graphs and erasures}.
Linear Algebra Appl. 404 (2005), 118--146.
\bibitem[BPS03]{BPS03}
Binz, Ernst; Pods, Sonja; Schempp, Walter, {\em Heisenberg groups---a unifying
structure of signal theory, holography and quantum information theory}. J. Appl.
Math. Comput. 11 (2003), no. 1-2, 1--57.


\bibitem[Dau92]{Dau92}
I. Daubechies, {\em Ten lectures on wavelets}. CBMS-NSF Regional Conf. Ser. in App. Math. 61. Soc. Industrial App. Math. (SIAM), Philadelphia, PA, 1992.
\bibitem[DuJo03]{DuJo03}
D.E. Dutkay, P.E.T. Jorgensen, {\em Wavelets on Fractals}; to appear
in Rev. Mat. Iberoamericana, preprint 2003;
http://arxiv.org/abs/math.CA/0305443
\bibitem[DuJo05]{DuJo05}
D.E. Dutkay, P.E.T. Jorgensen {\em Iterated function systems, Ruelle operators, and invariant projective measures}, preprint 2005, arxiv math.DS/0501077
\bibitem[Edg98]{Edg98}
G.A. Edgar, {\em Integral, probability, and fractal measures}, Springer-Verlag, New York, 1998
\bibitem[Er39]{Er39}
P. Erd\"os, {\em On a family of symmetric Bernoulli convolutions}.
Amer. J. Math. 61, (1939). 974--976.

\bibitem[Fal93]{Fal93}
K.J. Falconer, {\em Wavelets, fractals and order-two densities. Wavelets, fractals, and Fourier transforms}, 39--46, Inst. Math. Appl. Conf. Ser. New Ser., 43, Oxford Univ. Press, New York, 1993.

\bibitem[Ga62]{Ga62}
A.M. Garsia, {\em Arithmetic properties of Bernoulli convolutions}; Trans. Amer Math Soc 102 (1962) 409-432.

\bibitem[Gun00]{Gun00}
R. Gundy, {\em Low-pass filters, Martingales, and multiresolution
analysis}, Appl. Comput. Harmonic Anal. 9 (2000) 204-219.

\bibitem[Hun84]{Hun84}
T.W. Hungerford, {\em Algebra}, Springer-Verlag, NY 1984, first ed.
1974
\bibitem[Hut81]{Hut81}
J.E. Hutchinson, {\em Fractals and selfsimilarity}, Indiana Univ. Math.
J., 30 (1981), 713-747.

\bibitem[JeWi35]{JeWi35}
B. Jessen, A. Wintner, {\em Distribution functions and the Riemann zeta function}; Trans. Amer. Math. Soc. 38 (1935) 48--88.

\bibitem[JoPe96]{JoPe96}
P.E.T. Jorgensen, S. Pedersen, {\em Harmonic analysis of fractal
measures}, Constructive Approximation, v 12 (1996), 1-30.

\bibitem[JoPe98]{JoPe98}
P.E.T. Jorgensen, S. Pedersen, {\em Dense analytic subspaces in
fractal $L^2$-spaces}, J. d'Analyse Math 75, 1998, 185-228.

\bibitem[Jor04]{Jor04}
P.E.T. Jorgensen, {\em Iterated function systems, representations, and Hilbert space}, International J. Math.
15 (2004) 813-832.
\bibitem[Jor06]{Jor06}
Palle E. T. Jorgensen, {\em Analysis and Probability: Wavelets, signals,
fractals}, Springer-Verlag, Graduate Texts in Math, GTM vol 234, to appear in
Feb 2006.


\bibitem[Kah71]{Kah71}
J.-P.  Kahane, {\em Sur la distribution de certaines series aleatoires}; Colloque de Théorie des Nombres (Univ. Bordeaux, Bordeaux, 1969), 119--122. Bull. Soc. Math. France, Mém. No. 25, Soc. Math. France Paris, 1971.

\bibitem[Kat04]{Kat04}
Y. Katznelson, {\em An introduction to harmonic analysis}. Third edition. Cambridge Mathematical Library. Cambridge University Press, Cambridge, 2004.
\bibitem[KeWi35]{KeWi35}
R. Kershner, A. Wintner; {\em On symmetric Bernoulli convolutions}; Amer. J. Math. (1935),57, 541-548.


\bibitem[KSS03]{KSS03}
M. Keane, K. Simon, B. Solomyak, {\em The dimension of graph directed attractors with overlaps on the line
with an application to a problem in fractal image recognition}. Fund. Math. 180 (2003),  279--292.



\bibitem[LMW96]{LMW96}
K.-S. Lau, M.F. Ma, J. Wang,  {\em On some sharp regularity
estimates of $L^2$-scaling functions}, SIAM J. Math. Anal. 27
(1996), 835-864.




\bibitem[Nu01]{Nu01}
R.D. Nussbaum, {\em Periodic points of positive linear operators and Perron-Frobenius operators}; Int. Eq. and
Op. Th. 39 (2001) 41-97.
\bibitem[PeSo96]{PeSo96}
Y. Peres; B. Solomyak,  {\em Absolute continuity of Bernoulli convolutions, a simple proof}. Math. Res. Lett.
3 (1996), 231--239.

\bibitem[PeSo00]{PeSo00}
Y. Peres, B. Solomyak, {\em Existence of $L\sp q$ dimensions and entropy dimension for self-conformal
measures.} Indiana Univ. Math. J. 49 (2000), 1603--1621.


\bibitem[PSS00]{PSS00}
Y. Peres; W. Schlag; B. Solomyak, {\em Sixty years of Bernoulli convolutions; in Fractal geometry and stochastics, II} (Greifswald/Koserow, 1998), 39--65, Progr. Probab., 46, Birkhäuser, Basel, 2000.

\bibitem[Rud87]{Rud87}
W. Rudin, {\em Real and complex analysis}. Third edition. McGraw-Hill Book Co., New York, 1987. xiv+416 pp.

\bibitem[So95]{So95}
B. Solomyak, {\em On the random series $\sum\pm\lambda\sp n$ (an Erdos problem)}. Ann. of Math. (2) 142 (1995),
611--625.

\bibitem[SSU01]{SSU01}
K. Simon, B. Solomyak, M. Urbanski, {\em Invariant measures for parabolic IFS with overlaps and random
continued fractions.} Trans. Amer. Math. Soc. 353 (2001),  5145--5164.


\bibitem[Wal78]{Wal78}
P. Walters, {\em Equilibrium states for $\beta $-transformations and related transformations}. Math. Z. 159 (1978), no. 1, 65--88.
\bibitem[Wal01]{Wal01}
P. Walters, {\em Convergence of the Ruelle operator for a function satisfying Bowen's condition}. Trans. Amer. Math. Soc. 353 (2001), no. 1, 327-347.

\end{thebibliography}
\end{document}